\documentclass[reqno]{amsart}
\usepackage{hyperref}

\usepackage{a4wide,calrsfs}
\usepackage{amsmath,bbm}
\usepackage{color}
\usepackage{epsfig}
\usepackage{subfigure}
\usepackage{amssymb}
\usepackage{enumerate}
\usepackage{graphicx}
\usepackage{subfigure}
\usepackage{amsmath}
\usepackage{rotating}
\usepackage{stmaryrd}
\usepackage{upgreek}
\usepackage{cite}

\newcommand{\R}{\mathbb{R}}

\newcommand{\re}{\mathrm{Re}}

\allowdisplaybreaks

\begin{document}
\title[Nonoscillatory Behavior of Poincar\'e type equation]
{New results for the Nonoscillatory Asymptotic Behavior
of  High Order
Differential Equations of Poincar\'e type}

\author[A. Coronel \& F. Huancas]
{An{\'\i}bal Coronel and Fernando Huancas}  

\address{An{\'\i}bal Coronel, \newline
GMA, Departamento de Ciencias B\'asicas,\newline
Facultad de Ciencias, Universidad del B\'{\i}o-B\'{\i}o,\newline
Campus Fernando May, Chill\'{a}n, Chile}
\email{acoronel@ubiobio.cl}

\address{Fernando Huancas, \newline
GMA, Departamento de Ciencias B\'asicas,\newline
Facultad de Ciencias, Universidad del B\'{\i}o-B\'{\i}o,\newline
Campus Fernando May, Chill\'{a}n, Chile}
\email{fihuanca@gmail.com}

\thanks{\today}

\begin{abstract}
In this paper we study the asymptotic behavior 
of nonoscillatory solutions for high order
differential equations of Poincar\'e type.
We introduce two new and more weak than  classical hypotheses on the coefficients,
which implies a well posedness result and
a characterization of asymptotic behavior for the solution
of the Poincar\'e equation. 
The proof of the main results is based on the application of the scalar method.
We define a change 
of variable to reduce the order of Poincar\'e type equation
and deduce that the new variable
satisfies a nonlinear differential equation.
We apply the variation of parameters method and the Banach fixed point 
theorem to get the well posedness and asymptotic behavior 
of the nonlinear equation. Then, by rewritten the results in terms
of the original variable we establish the existence 
of a fundamental system of solutions for and we precise 
the formulas for the asymptotic behavior of the Poincar\'e type equation.

%

\end{abstract}

\maketitle
\numberwithin{equation}{section}
\newtheorem{theorem}{Theorem}[section]
\newtheorem{corollary}{Corollary}[theorem]
\newtheorem{proposition}[theorem]{Proposition}
\newtheorem{definition}{Definition}[section]
\newtheorem{lemma}{Lemma}[section]
\allowdisplaybreaks

\section{Introduction}

\subsection{Scope}
\label{sec:scope}
In this paper we are interested in the nonoscillatory  asymptotic behavior of 
the following differential equation:
\begin{align}
y^{(n)}+\sum_{i=0}^{n-1}[a_i+r_i (t)] y^{(i)}
&=0,
\quad n\in \mathbb{N}, n\ge 2,
\label{eq:intro_uno_intro}
\end{align}
where $a_i\in\R$ are given constants  and 
$r_i$ are real-valued functions. 
The equation 
\eqref{eq:intro_uno_intro} is called 
of a Poincar\'e type, since its analysis was
motivated and introduced by Poincar\'e \cite{poincare_1885}.
Then, naturally, the question of existence, uniqueness and  asymptotic behavior
of \eqref{eq:intro_uno_intro} is called the Poincar\'e problem. 
Later on the Poincar\'e work the asymptotic behavior of 
\eqref{eq:intro_uno_intro} have been investigated  by several authors with a
long and rich history of results 
\cite{bellman_book,eastham_book,harris_lutz_1977,hartman_1948}.
For a recent brief review
of some historical landmarks of the evolution of Poincar\'e problem
we refer to \cite{coronel_huancas_pinto}. 

The purpose of this paper is to solve the Poincar\'e problem by introducing more weak
than classical hypotheses for perturbation functions $r_i$. Our analysis is based on the
modified scalar method which consists in three big steps, for
details on cases $n=3,4$ see  \cite{figueroa_2008} and \cite{coronel_huancas_pinto},
respectively. In the first step we introduce the change of variable
\begin{align}
z(t)=\frac{y'(t)}{y(t)}-\mu 
\quad\mbox{or equivalently}\quad
y(t)=\exp\Big(\int_{t_0}^t (z(s)+\mu )ds\Big).
\label{eq:general_change_var}
\end{align}
to reduce the order of \eqref{eq:intro_uno_intro}.
If $\mu\in\R$ is a simple characteristic
root of \eqref{eq:intro_uno_intro} when $r_i=0$ then 
$z$ satisfy a nonlinear differential equation of the following
type
\begin{align}
z^{(n-1)}(t)+ \sum_{i=0}^{n-2} b_{i}(\mu)z^{(i)}(t)
&=\mathbb{P}(\mu,t,r_0(t),\ldots,r_{n-1}(t),z(t),z'(t),\ldots,z^{(n-2)}(t)),
\label{eq:riccati_type_general}
\end{align}
where $b_i$ are real-valued functions,
 $\mu\in\R$ is a given (fix)  parameter,  
and $\mathbb{P}:\R^{2n+1}\to\R$ is a polynomial of $n$ degree in the 
$n-1$ last variables and the coefficients depends on
the first $n+2$ variables, i.e.
 $\mathbb{P}$ admits the representation 
\begin{subequations}
 \label{eq:general_no_lin:2}
 \begin{align}
&\mathbb{P}(\mathbf{e},\mathbf{x})=\sum_{|\boldsymbol{\alpha}|=0}^{n}
\Omega_{\boldsymbol{\alpha}}(\mathbf{e})\mathbf{x}^{\boldsymbol{\alpha}}
\quad
\mbox{with
$\mathbf{e}=(e_1,\ldots,e_{n+2})$
and 
$\mathbf{x}=(x_1,\ldots,x_{n-1})$},
\label{eq:general_no_lin:2_1}
\\
&\boldsymbol{\alpha}=(\alpha_1,\ldots,\alpha_{n-1})\in\mathbb{N}^{{n-1}}_0,
\quad \mathbb{N}_0=\mathbb{N}\cup\{0\}.
\label{eq:general_no_lin:2_2}
\end{align}
\end{subequations}
Here, we have used the standard multindex notation for 
$|\boldsymbol{\alpha}|$ and $\mathbf{x}^{\boldsymbol{\alpha}}$, i.e. 
$|\boldsymbol{\alpha}|=\alpha_1+\ldots+\alpha_{n-1}$ and 
$\mathbf{x}^{\boldsymbol{\alpha}}=x_1^{\alpha_1}\ldots x_3^{\alpha_{n-1}}$.
The equation \eqref{eq:riccati_type_general} is the abstract form 
of a more particular equation obtained in the specific case of change of
variable \eqref{eq:general_change_var}, for instance in the case $n=3,4$
there is several coefficients such that 
$\Omega_{\boldsymbol{\alpha}}=0$ \cite{figueroa_2008,coronel_huancas_pinto}.
In the second step we analyze  the existence, uniqueness and asymptotic behavior
of \eqref{eq:riccati_type_general}. Finally, in a third step we particularize
the results of \eqref{eq:riccati_type_general} to the specific equation
obtained for $\mu$ characteristic root of \eqref{eq:intro_uno_intro}
and using the definition of $y$ in \eqref{eq:riccati_type_general}
we deduce the well posedness and asymptotic behavior of \eqref{eq:intro_uno_intro}.

The main aims of the paper are the following
\begin{itemize}
 \item[(O$_1$)] Provide a well posedness framework to analyze \eqref{eq:riccati_type_general}.
 \item[(O$_2$)] Precise the asymptotic behavior of \eqref{eq:riccati_type_general}.
 \item[(O$_3$)] Introduce the appropriate assumptions such that the Poincar\'e problem
 can be solved by particularize the results of (O$_1$) and (O$_2$).
\end{itemize}
Now, in order to precise the results for each objective we 
introduce some notation in subsection 
\ref{sec:notPOI} and \ref{sec:notRIC} 
and present the results on subsection \ref{sec:main_results}.

\subsection{Assumptions for the coefficients of \eqref{eq:riccati_type_general}}
\label{sec:notPOI}
Let us introduce some notation. 
Given  the Green function $g$ (see \eqref{eq:green_function}),
we consider that the  
functions $\mathcal{R},\mathcal{L}_k$ and $\Upsilon_{\ell}$
 are defined as follows
\begin{align}
&\mathcal{R}(t)=
\sum_{j=0}^{n-2}
\left|\int_{t_0}^{\infty}\frac{\partial^j g}{\partial t^j}(t,s)
	\Omega_{\boldsymbol{0}}\Big(\mu,s,r_0(s),\ldots,r_{n-1}(s)\Big)ds\right|,
\label{eq:operator_l0}
\\
&\mathcal{L}_k(t)=\int_{t_0}^{\infty}\left[\sum_{j=0}^{n-2}
	\left|\frac{\partial^j g}{\partial t^j}(t,s)\right|
	\right]\sum_{|\boldsymbol{\alpha}|=k}\Big|
	\Omega_{\boldsymbol{\alpha}}
	\Big(\mu,s,r_0(s),\ldots,r_{n-1}(s)\Big)\Big|ds,
	k=\overline{1,n}
\label{eq:operator_lk}
\\
&\Upsilon_0(y_1,\ldots,y_d)=
\prod_{1\le i<j\le d}(y_j-y_i),
\quad 
\Upsilon_{\ell}(y_1,\ldots,y_d)
=\!\!\!\!
\prod_{
\begin{array}{c}
1\le i<j\le d \\
i\not=\ell,j\not=\ell
\end{array}
}
\!\!\!\!
(y_j-y_i)
\;\;\;\;
\mbox{for }
\ell=\overline{1,d}.
\label{eq:phi_notation}
\end{align}
Then, we assume the following hypotheses on $b_i(\mu)$ and $\mathbb{P}$:
\begin{enumerate}
\item[(R1)]  The functions $b_i(\mu)$ for $i=\overline{0,n-2}$ are such that the 
 set of characteristic roots for \eqref{eq:riccati_type_general} when $\mathbb{P}=0$
is given by $\Gamma_\mu=\Big\{\gamma_i(\mu), i=\overline{1,n-1}\;:\;
\gamma_1>\gamma_2>\ldots>\gamma_{n-1}\Big\}\subset\R.$
 \item [(R2)] The functions $\Omega_{\boldsymbol{\alpha}}$ given on
\eqref{eq:general_no_lin:2_1}
satisfy the following requirements
\begin{align}
 \lim_{t\to\infty}\mathcal{R}(t)
 =\lim_{t\to\infty}\mathcal{L}_1(t)=0
 \quad\mbox{and}\quad 
 \lim_{t\to\infty}\sum_{k=2}^n\mathcal{L}_k(t)<1,
 \label{eq:integral_codition}
\end{align}
where $\mathcal{R}$ and $\mathcal{L}_k$ are defined on \eqref{eq:operator_lk}.
 \item [(R3)] 
 The functions $\Omega_{\boldsymbol{\alpha}}$  on
\eqref{eq:general_no_lin:2_1}
satisfy the following property: given $\Phi_{1}$ and $\gamma_i$
with
\begin{align}
\Phi_{1}=\Big[\Big|\Upsilon_0(\gamma_1,\ldots,\gamma_{n-1})\Big|\Big]^{-1}
\sum_{\ell=1}^{n-1}
\Big|\Upsilon_{\ell}(\gamma_1,\ldots,\gamma_{n-1})\Big|
\sum_{j=0}^{n-2}(|\gamma_{\ell}|)^j,
\label{eq:PHI_1_1}
\end{align}
and $\gamma_i\in \Gamma_{\mu}$, then there exits 
$\sigma_{\gamma_i}\in [0,1/\Phi_1[$ such that the inequality
\begin{align*}
 \int_{t_0}^{\infty} \exp(-\gamma_i(t-s))
\sum_{|\boldsymbol{\alpha}|=1}^n
\Big|\Omega_{\boldsymbol{\alpha}}\Big(\mu,s,r_0(s),\ldots,r_{n-1}(s)\Big)\Big|
ds\le \sigma_{\gamma_i}
\end{align*}
is satisfied for all $t\in [t_0,\infty[$.
\end{enumerate}
The assumptions (R1) and (R2) are considered 
to study the existence and the uniqueness
of the solution for \eqref{eq:riccati_type_general}. The hypothesis (R1)
is considered, since we have interested in the nonoscillatory behavior.
Meanwhile, (R2) is the natural condition in order to get the application of Banach
fixed point theorem. Now,  the hypothesis (R3) is used in order to get  the asymptotic behavior 
of the solution for \eqref{eq:riccati_type_general}.

\subsection{Assumptions for the coefficients of \eqref{eq:intro_uno_intro}}
\label{sec:notRIC}
Let us consider $\gamma_j(\lambda_i)$ and the function $\mathbb{H}$ defined by
 \begin{align}
& \gamma_j(\lambda_i)=
\left\{
\begin{array}{ll}
 \lambda_j-\lambda_i,& j=1,\ldots,i-1,
 \\
 \lambda_{j+1}-\lambda_i,& j=i,\ldots,n-1,
\end{array}
\right.
\label{eq:espectral} 
\\
&\mathbb{H}(t,\mu)=(n-1)+a_2+r_2(t)\Big(2+2\mu\Big)+r_1(t)+\sum_{m=1}^{n-3}\hat{S}_{m,n}
\nonumber\\
&
\qquad
+\sum_{i=3}^{n-1} \left[a_i\Big((i-1)
+\sum_{m=1}^{i-3}\hat{S}_{m,i}+1+\sum_{j=1}^{i-2}{i\choose j}
\mu^{j}\Big) 
+
r_i(t)\Big(\sum_{m=0}^{i-3}\hat{S}_{m,i}(t)+(1+\mu)^{i}\Big)
\right]
\label{eq:function:HP},
\end{align}
where $\hat{S}_{0,j}(t)=j+\mu (j-1)$ and for $m> 1$ we have that
\begin{align*}
 \hat{S}_{m,j}&=\sum_{\ell_1=0}^{j-(m+2)}
\quad\sum_{\ell_2=0}^{j-\ell_1-(m+2)}
\ldots
\sum_{\ell_{m}=0}^{j-\sum_{q=0}^{m-1} \ell_q-(m+2)}
{j-1\choose \ell_{1}}
\prod_{i=1}^{m-1}{j-\sum_{q=1}^i \ell_q-i-1\choose \ell_{i+1}}
\nonumber\\
&
\hspace{2cm}
\times
\Bigg[1
+\Big(j-\sum_{q=1}^m \ell_q-m-1\Big)
\Big(1+\mu\Big)\Bigg].
\end{align*}
Then,
related with the coefficients of \eqref{eq:intro_uno_intro}, we consider that:
\begin{enumerate}
 \item[(H1)] The constants $a_i$ are selected such that
 the set of characteristic  roots for \eqref{eq:intro_uno_intro}
 when $r_i=0$ is given 
by $\Lambda=\Big\{\lambda_i, i=\overline{1,n}\;:\;
\lambda_1>\lambda_2>\ldots>\lambda_n\Big\}\subset\R.$
 \item[(H2)] The perturbation functions 
 satisfy the asymptotic behavior:
\begin{align}
 &\displaystyle 
 \lim_{t\to\infty}\sum_{i=0}^{n-2} \left|\int_{t_0}^{\infty}
  \frac{\partial^i g}{\partial t^i}(t,s)r_j(s)ds\right|=0
  \quad\mbox{for $j=\overline{0,n-1}$},
  \\
 &\displaystyle
  \lim_{t\to\infty}\int_{t_0}^{\infty}\sum_{i=0}^{n-2}
  \left|\frac{\partial^i g}{\partial t^i}(t,s)\right|
 \left|\sum_{\ell=0}^n\mu^{\ell}r_{\ell}(t)\right|ds =0
  \quad\mbox{for each  $\mu\in\Lambda$,}
\label{eq:H2}
\end{align}
for a Green function $g$.

 \item[(H3)] The perturbation functions $r_i$ and the characteristic set  
 $\Lambda$
satisfy the following property: given $\Phi_{i}^\Lambda$ and $\gamma_k(\lambda_i)$
with
\begin{align*}
\Phi_{i}^\Lambda=\Big[\Big|\Upsilon_0(\gamma_1(\lambda_i),\ldots,
  \gamma_{n-1}(\lambda_i))\Big|\Big]^{-1}
\sum_{\ell=1}^{n-1}
\Big|\Upsilon_{\ell}(\gamma_1(\lambda_i),\ldots,\gamma_{n-1}(\lambda_i))\Big|
\sum_{j=0}^{n-2}(|\gamma_{\ell}(\lambda_i)|)^j,
\end{align*}
and $\gamma_k(\lambda_i)\in \Gamma_{\lambda_i}$, then there exits 
$\sigma_{\gamma_k}\in [0,1/\Phi^\Lambda_i[$ such that the inequality
\begin{align*}
 \int_{t_0}^{\infty} \exp(-\gamma_k(\lambda_i)(t-s))
\Big|\mathbb{H}(s,\lambda_i)\Big|
ds\le \sigma_{\gamma_k}
\end{align*}
is satisfied for all $t\in [t_0,\infty[$.

\end{enumerate}

\subsection{Main results}
\label{sec:main_results}

The results related with the aims (O$_1$), (O$_2$) and (O$_3$)
(see subsection \ref{sec:scope})
are given below on Theorem~\ref{teo:ricati},  
Theorem~\ref{teo:asymptotic_ricatti}, and
Theorems~\ref{teo:changeofvariable}-\ref{teo:main}, respectively.

\begin{theorem}
\label{teo:ricati}
For $t_0\in\R$ (fix) 
consider the notation $C_0^{n-2}([t_0,\infty[)$
for set of functions
\begin{align*}
C_0^{n-2}([t_0,\infty[)
=\left\{z\in C^{n-2}([t_0,\infty[,\R)\quad :\quad
\lim_{t\to\infty}z^{(k)}(t)= 0\mbox{ for }
k=\overline{0,n-2}\right\},
\end{align*}
which is a Banach space with the norm $
\|z\|_0=\sup_{t\geq t_0}\sum_{i=0}^{n-2}|z^{(i)}(t)|$.
Assume that the coefficients of
the equation \eqref{eq:riccati_type_general}
satisfy (R1) and (R2).
Then, exists a unique  $z\in C_0^{n-2}([t_0,\infty[)$ solution 
of~\eqref{eq:riccati_type_general}.
\end{theorem}

\begin{theorem}
\label{teo:asymptotic_ricatti}
Let us introduce  the notation
\begin{align}
\mathbb{G}_\mu&=\Big\{(\gamma_1,\ldots,\gamma_{n-1})\;\;:\;\; 
\gamma_i\in \Gamma_\mu, i=\overline{1,n-1}\Big\},
\label{eq:notation:G}
\\
 \mathbb{E}^{n-1}_i&=
\Big\{(x_1,\ldots,x_{n-1})\in \R^{n-1}\;:\; 
x_1>x_2>\ldots >x_{n-1}\quad\mbox{and}\quad
\nonumber\\
&\hspace{1cm}
\Big(\mbox{$x_1<0$ if $i=1$ or }
 \mbox{$0\in ]x_i,x_{i-1}[$ if $i=\overline{2,n-1}$ or  $x_{n-1}>0$
if $i=n$}\Big)
\Big\}.
\label{eq:notation:E}
\end{align}
Consider that the hypotheses of Theorem~\ref{teo:ricati}
and the assumption (R3) are valid. Then,  $z_{\mu}$ the solution
of \eqref{eq:riccati_type_general} has 
the following asymptotic behavior 
\begin{align}
z^{(j)}_{\mu}(t)
=
\left\{
\begin{array}{l}
\displaystyle
O\left(\int_{t}^{\infty}e^{-\beta(t-s)}
\Big|\Omega_{\boldsymbol{0}}
\Big(\mu,s,r_0(s),\ldots,r_{n-1}(s)\Big)\Big|ds\right),
\\
\hspace{3cm} \mbox{when }\mathbb{G}_\mu\subset\mathbb{E}^{n-1}_1\mbox{ and } \beta\in [\gamma_1,0[,
\\
\displaystyle
O\left(\int_{t_0}^{\infty}e^{-\beta(t-s)}
\Big|\Omega_{\boldsymbol{0}}\Big(\mu,s,r_0(s),
\ldots,r_{n-1}(s)\Big)\Big|ds\right),
\quad
\\
\hspace{3cm} 
\mbox{when }
\mathbb{G}_\mu\subset\mathbb{E}^{n-1}_i \mbox{ and } \beta\in [\gamma_i,0[,\quad 
 i=\overline{2,n-1},
\\
\displaystyle
O\left(\int_{t_0}^{t}e^{-\beta(t-s)}
\Big|\Omega_{\boldsymbol{0}}\Big(\mu,s,r_0(s),\ldots,r_{n-1}(s)\Big)\Big|ds
\right),
\\
\hspace{3cm}
\mbox{when }
\mathbb{G}_\mu\subset\mathbb{E}^{n-1}_n\mbox{ and } \beta\in ]0,\gamma_{n-1}],
\end{array}
\right.
\label{eq:asymptotic_teo:form}
\end{align}
when $t\to\infty,$
for all $j\in\{0,\ldots,n-2\}$.
\end{theorem}

\begin{theorem}
\label{teo:changeofvariable} 
Let us consider the new variable $z$ satisfying
\eqref{eq:general_change_var}.
If $y$ is a solution of \eqref{eq:intro_uno_intro},
the following assertions are valid:

\vspace{0.3cm}
\noindent (a)
If $\mu\in\Lambda$, then
$z$ satisfy the differential equation 
\begin{align}
&z^{(n-1)}+\sum_{i=2}^{n-1}\Big[\sum_{k=i}^{n}a_k {k \choose i} \mu ^{k-i} \Big]z ^{(i-1)}
+\Big(n\mu ^{n-1}+\sum_{i=1}^{n-1}i a_i \mu ^{i-1}\Big) z
=-\mathbb{F}
 \label{eq:ricati_enesima}
\end{align}
where
\begin{align}
\mathbb{F}&=(n-1)z^{(n-2)}(t)z(t)+\sum_{m=1}^{n-3}S_{m,n}(t)
\nonumber\\
&\quad
+
\sum_{i=3}^{n-1} \left[a_i\Big((i-1)z^{(i-2)}(t)z(t)
+\sum_{m=1}^{i-3}S_{m,i}(t)+z^{(i-1)}(t)z^{(1)}(t)+\sum_{j=1}^{i-2}{i\choose j}
z^{(i-j)}(t)\mu^{j}\Big) 
\right.
\nonumber\\
&\hspace{1.8cm}
\left.
+
r_i(t)\Big(\sum_{m=0}^{i-3}S_{m,i}(t)+z^{(i-1)}(t)z^{(1)}(t)+(z(t)+\mu)^{i}\Big)
\right]
\nonumber\\
&\quad
+a_2z^2(t)+r_2(t)\Big(z^{(1)}(t)+(z(t)+\mu)^2\Big)
+r_1(t)\Big(z(t)+\mu\Big)+r_0(t),
\label{eq:function:FF}\\
 S_{0,j}(t)&=z^{(j-1)}(t) +(j-1)\;z^{(j-2)}(t)\; (z(t)+\mu)
 \label{eq:leibiniz:j_0}\\
 S_{m,j}(t)&=\sum_{\ell_1=0}^{j-(m+2)}
\quad\sum_{\ell_2=0}^{j-\ell_1-(m+2)}
\ldots
\sum_{\ell_{m}=0}^{j-\sum_{q=0}^{m-1} \ell_q-(m+2)}
{j-1\choose \ell_{1}}
\prod_{i=1}^{m-1}{j-\sum_{q=1}^i \ell_q-i-1\choose \ell_{i+1}}
\nonumber\\
&
\quad
\times
\prod_{i=1}^{m}(z+\mu)^{(\ell_{i})}(t)
\Bigg[(z+\mu)^{(j-\sum_{q=1}^m \ell_q-m-1)} (t)
\nonumber\\
&
\hspace{3cm}
+\Big(j-\sum_{q=1}^m \ell_q-m-1\Big)
(z+\mu)^{(j-\sum_{q=1}^m \ell_q-m-2)} (t)
\Big(z(t)+\mu\Big)\Bigg].
\label{eq:leibiniz:j_m}
\end{align}

\vspace{0.3cm}
\noindent (b) If the hypothesis (H1) is satisfied, then
\begin{align}
\Gamma_{\lambda_i}=\Big\{\gamma_j\;:\;
\gamma_j(\lambda_i)\mbox{ given on \eqref{eq:espectral}, }\;\;
\gamma_1>\gamma_2>\ldots>\gamma_{n-1}\Big\}\subset\R
\label{eq:charset_ricca}
\end{align}
is the characteristic set of \eqref{eq:ricati_enesima}
when $\mu=\lambda_i$.

\vspace{0.3cm}
\noindent (c) If the hypothesis (H1) and (H2)
are satisfied and $\mu=\lambda_i$. Then, exists a unique 
$z_{\lambda_i}\in C_0^{n-2}([t_0,\infty[)$ solution 
of \eqref{eq:ricati_enesima}.

\vspace{0.3cm}
\noindent (d) If the hypothesis (H1), (H2) and (H3)
are satisfied. Then
$z_{\lambda_i}$, the solution of \eqref{eq:ricati_enesima}
with $\mu=\lambda_i$,
satisfies  the following asymptotic behavior
\begin{align}
z^{(j)}_{\lambda_i}(t)
=
\left\{
\begin{array}{l}
\displaystyle
O\left(\int_{t}^{\infty}e^{-\beta(t-s)}
\Big|\sum_{\ell=0}^{n-1}\lambda^{\ell}_1r_{\ell}(s)\Big|ds\right),
\;\; i=1,\;\;\beta\in [\lambda_2-\lambda_1,0[,
\\
\displaystyle
O\left(\int_{t_0}^{\infty}e^{-\beta(t-s)}
\Big|\sum_{\ell=0}^{n-1}\lambda^{\ell}_ir_{\ell}(s)\Big|ds\right),
i=\overline{2,n-1},\;\;\beta\in [\lambda_{i+1}-\lambda_i,0[,
\\
\displaystyle
O\left(\int_{t_0}^{t}e^{-\beta(t-s)}
\Big|\sum_{\ell=0}^{n-1}\lambda^{\ell}_nr_{\ell}(s)\Big|ds
\right),
\;\; i=n,\;\; \beta\in ]0,\gamma_{n-1}],
\end{array}
\right.
\label{eq:asymptotic_teo:form:lambda_i}
\end{align}
when $t\to\infty,$
for all $j\in\{0,\ldots,n-2\}$.

\end{theorem}

\begin{theorem}
\label{teo:main}
Let us consider that (H1) and   (H2) are satisfied.
Then, the equation \eqref{eq:intro_uno_intro} has a fundamental system of 
solutions given
\begin{align}
    y_i(t)=\exp\Big(\int_{t_0}^{t}[\lambda_i+z_{\lambda_i}(s)]ds\Big),
    \quad
    \mbox{with $z_{\lambda_i}$ solution of \eqref{eq:ricati_enesima}
    with $\mu=\lambda_i\in\Lambda$}.
    \label{eq:fundam_sist_sol}
\end{align}
 Moreover, the following properties about
 the asymptotic behavior
\begin{align}
\lim_{t\to\infty}
\frac{y^{(j)}_i(t)}{y_i(t)}&=(\lambda_i)^j\quad j=\overline{1,n},
\label{eq:asymptotic_perturbed_w0}\\
W[y_1,\ldots,y_n](t)&= 
\prod_{1\le k<\ell\le n}\;\big(\lambda_{\ell}-\lambda_k\big)
\; \prod_{i=1}^n y_i(t)(1+o(1)),
\label{eq:asymptotic_perturbed_w}
\end{align}
is satisfied when $t\to\infty.$
Furthermore, if $\pi_i=\prod_{1\le i<j\le n, j\not=i}(\lambda_j-\lambda_i)$
and (H3) is satisfied, then 
\begin{align}
&y^{(j)}_i(t)=
\Big(\lambda^j_i+o(1)\Big)
e^{\lambda_i(t-t_0)}
\nonumber\\
&\hspace{1.5cm}
\times
	\exp\left(
	\frac{1}{\pi_i}
	\int_{t_0}^{t}
	\mathbb{F}\Big(\lambda_i,s,r_0(s),\ldots,r_{n-1}(s),
z_{\lambda_i}(s),\ldots,z^{(n-2)}_{\lambda_i}(s)\Big)ds\right),
\label{eq:asymptotic_perturbed_3_4}
\end{align}
holds, when $t\to\infty$ with 
$z^{(j)}_{\lambda_i},\;\;j=\overline{0,n-2}$ 
given asymptotically by \eqref{eq:asymptotic_teo:form:lambda_i} and $\mathbb{F}$ 
defined on \eqref{eq:function:FF}.
\end{theorem}

\subsection{A brief discussion of the assumptions.}
Related with the constant coefficients part of \eqref{eq:intro_uno_intro}, 
 i. e. when $r_i=0$, there is a coincidence or 
 common hypothesis used by the different researchers, in order
to deduce the nonoscillatory asymptotic behavior. All of them consider 
the fact that 
\begin{align}
\left.
\begin{array}{l}
\mbox{There is a simple characteristic root $\mu$
such that}
\\\mbox{$\re(\mu)\not=\re(\mu_0)$
for any other characteristic root $\mu_0$.}
\end{array}
\right\}
\label{H0:old}
\end{align}
Then, (H1) is an 
extension of the condition \eqref{H0:old} for all characteristic roots
of \eqref{eq:intro_uno_intro}.

On the other hand,
with respect to the regularity and the asymptotic behavior of perturbation functions $r_i$
we have different assumptions. 
For instance, the seminal work of Poincar\'e \cite{poincare_1885}
consider that $r_i$ are rational functions,  
Perron \cite{perron_1909} assumes that $r_i$ are continuous
functions, Levinson \cite{levinson1948} considers that
$r_i\in L^1([t_0,\infty[)$, and 
Hartman and Wintner in \cite{hartman_wintner_1955}
select $r_i\in L^p([t_0,\infty[)$ for some $p\in ]1,2]$.
There is a coincidence 
of the authors with respect to the asymptotic behavior of $r_i$,
by considering that:  $r_i(t)\to 0$ when $t\to\infty$.
Then (H2) and (H3) are new and more weak than  the previous assumptions.
Indeed, we refer to the recent work \cite{coronel_huancas_pinto}, in the 
case $n=4$,
for an  example of perturbation functions $r_i$
which satisfy (H2)-(H3) and neither satisfy the
assumptions of $L^p$ regularity and nor
satisfy the asymptotic behavior $r_i(t)\to 0$ when $t\to\infty$.

\section{Proof of main results}
\label{sec:proofs_mr}

\subsection{Proof of Theorem~\ref{teo:ricati}}
The proof is mainly based on variation of parameters technique and Banach fixed
point Theorem.  Indeed, we proceed as follows:  we introduce the 
precise notation of
Green functions  for \eqref{eq:riccati_type_general}
when $\mathbb{P}=0$,  we apply the method of variation 
of parameters to get the operator equation, and  
we deduce that the operator satisfies the hypotheses of Banach fixed
point Theorem.

The Green function $g$ for \eqref{eq:riccati_type_general}
when $\mathbb{P}=0$ is defined by 
\begin{align}
g(t,s)=\frac{-g_{\mu}(t,s)}{\Upsilon_0(\gamma_1,\ldots,\gamma_{n-1})}
\quad
\mbox{when
 $\mathbb{G}_{\mu}\subset\mathbb{E}^{n-1}_i$
 for $i=\overline{1,n}$,}
\label{eq:green_function}
\end{align} 
where $\Upsilon_0$ is the notation on \eqref{eq:phi_notation},
$\mathbb{G}_{\mu}$ is the set defined on \eqref{eq:notation:G}, the 
notation $\mathbb{E}^{n-1}_i$ is given on \eqref{eq:notation:E},
and $g_{\mu}:\R^2\to\R$
are the  functions defined as follows
\begin{align*}
 g_{\mu}(t,s)=&
\left\{
\begin{array}{lcl}
\displaystyle
\sum_{\ell=1}^{n-1} G_{\ell}e^{-\gamma_\ell(t-s)}H(s-t), &  &\mathbb{G}_{\mu}\subset\mathbb{E}^{n-1}_1,
\\
\displaystyle
\sum_{\ell=1}^{k-1} G_{\ell}e^{-\gamma_\ell(t-s)} H(t-s)
+\sum_{\ell=k}^{n-1} G_{\ell}e^{-\gamma_\ell(t-s)}H(s-t), &  &
\mathbb{G}_{\mu}\subset\mathbb{E}^{n-1}_k,\,\,\,k=\overline{2,n-1},
\\
\displaystyle
\sum_{\ell=1}^{n-1} G_{\ell}e^{-\gamma_\ell(t-s)}H(t-s), &  &
\mathbb{G}_{\mu}\subset\mathbb{E}^{n-1}_n,
\end{array}
\right.
\end{align*}
where $G_\ell=(-1)^\ell\Upsilon_{\ell}(\gamma_1,\ldots,\gamma_{n-1})$ 
with $\Upsilon_\ell$ defined on \eqref{eq:phi_notation}
and
$H$ is the Heaviside function, i.e. $H(x)=1$ for $x\ge 0$ and 
$H(x)=0$ for $x<0.$ For further details on Green functions the reader
may be consult~\cite{bellman_book}
(see also \cite{figueroa_2006} for $n=2$,
\cite{figueroa_2008} for $n=3$ and \cite{coronel_huancas_pinto} for $n=4$).

We apply the method of variation 
of parameters to get  that
\eqref{eq:riccati_type_general} is equivalent to  the 
following integral equation
\begin{align}
z(t)=\int_{t_0}^\infty g(t,s)\mathbb{P}\Big(\mu,s,r_0(s),
\ldots,r_{n-1}(s),z(s),\ldots,z^{(n-2)}(s)\Big)ds,
\label{eq:integ_equation}
\end{align}
where $g$ is the Green function defined on  \eqref{eq:green_function}.
Thus, if we define the operator $T$ from $C_0^{n-2}([t_0,\infty[)$ to 
$ C_0^{n-2}([t_0,\infty[)$ as follows
\begin{align}
Tz(t)&=\int_{t_0}^\infty g(t,s)\mathbb{P}\Big(\mu,s,r_0(s),
\ldots,r_{n-1}(s),z(s),\ldots,z^{(n-2)}(s)\Big)ds,
\label{eq:operator_fix_point}
\end{align}
 we note that \eqref{eq:integ_equation} can be rewritten as the operator
equation 
\begin{align}
Tz=z
\qquad
\mbox{over}
\qquad
D_{\upeta}:=\Big\{z\in C_0^{n-2}([t_0,\infty[)\quad:\quad \|z\|_0\le\upeta \Big\},
\label{eq:operator_equation}
\end{align}
where $\upeta\in \R^+$ will be selected
in order to apply the Banach fixed point theorem.
Indeed, we have that

\vskip 0.5cm
\noindent
{\bf (a)} {\it $T$ is well defined from $C_0^{n-2}([t_0,\infty[)$ to $C_0^{n-2}([t_0,\infty[)$}.
Let us consider an arbitrary $z\in C_0^{n-2}([t_0,\infty[),$ 
by the definition of the operator $T$ we deduce that
 \begin{align*}
  T^{(j)}z(t) &= 
  \int_{t_0}^{\infty}
	\frac{\partial^j g}{\partial t^j}(t,s)
	\mathbb{P}\Big(\mu,s,r_0(s),\ldots,r_{n-1}(s),z(s),\ldots,z^{(n-2)}(s)\Big)ds,
	\quad
	j=\overline{0,n-2}.
\end{align*}
Then,  by the hypothesis (R2) and
using the notation \eqref{eq:general_no_lin:2},
we can deduce the following estimates
\begin{align}
&\sum_{j=0}^{n-2}|T^{(j)}z(t)| 
=
\sum_{j=0}^{n-2}
\left|\int_{t_0}^{\infty}
	\frac{\partial^j g}{\partial t^j}(t,s)
	\sum_{|\boldsymbol{\alpha}|=0}^{n}
	\Omega_{\boldsymbol{\alpha}}\Big(\mu,s,r_0(s),\ldots,r_{n-1}(s)\Big)
	\Big(z(s),\ldots,z^{(n-2)}(s)\Big)^{\boldsymbol{\alpha}}ds\right|
\nonumber\\
&\quad\le  \sum_{j=0}^{n-2} \left|\int_{t_0}^{\infty}
	\frac{\partial^j g}{\partial t^j}(t,s)
	\Omega_{\boldsymbol{0}}\Big(\mu,s,r_0(s),\ldots,r_{n-1}(s)\Big)ds\right|
\nonumber\\
&
\qquad
	+
	\int_{t_0}^{\infty}
	\sum_{j=0}^{n-2}
	\left|\frac{\partial^j g}{\partial t^j}(t,s)\right|
	\left|\sum_{|\boldsymbol{\alpha}|=1}^{n}
	\Omega_{\boldsymbol{\alpha}}\Big(\mu,s,r_0(s),\ldots,r_{n-1}(s)\Big)
	(z(s))^{\alpha_1}\ldots(z^{(n-2)}(s))^{\alpha_{n-1}}\right|ds
\nonumber\\
&
\quad\le \sum_{j=0}^{n-2} \left|\int_{t_0}^{\infty}
	\frac{\partial^j g}{\partial t^j}(t,s)
	\Omega_{\boldsymbol{0}}\Big(\mu,s,r_0(s),\ldots,r_{n-1}(s)\Big)ds\right|
\nonumber\\
&
\qquad
	+\sum_{k=1}^n
	\Big(\|z\|_0\Big)^k
	\int_{t_0}^{\infty}
	\sum_{j=0}^{n-2}
	\left|\frac{\partial^j g}{\partial t^j}(t,s)\right|
	\sum_{|\boldsymbol{\alpha}|=k}\Big|
	\Omega_{\boldsymbol{\alpha}}\Big(\mu,s,r_0(s),\ldots,r_{n-1}(s)\Big)\Big|ds
\label{eq:bound_for_tz0:prev}\\
&
\quad
\le\mathcal{R}(t)+\|z\|_0\mathcal{L}_1(t)
+\Big(\|z\|_0\Big)^2\max\Big\{1,\|z\|_0,\ldots,\Big(\|z\|_0\Big)^{n-2}\Big\}
\sum_{k=2}^n\mathcal{L}_k(t),
\label{eq:bound_for_tz0}
\end{align}
for $j=0,\ldots,n-2$.
Now, by \eqref{eq:integral_codition}
and the fact that  $z\in C^{n-2}_0([t_0,\infty[)$, we have that
the right hand side of 
\eqref{eq:bound_for_tz0} tend to $0$
when $t\to\infty$. Thus, for all $j=0,\ldots,n-2$
we have that $T^{(j)}z(t)\to 0$ when $t\to\infty$ 
or equivalently $Tz\in C^{n-2}_0([t_0,\infty[)$ 
for all $z\in C^{n-2}_0([t_0,\infty[).$

\vskip 0.5cm
\noindent
{\bf (b)} {\it For all $\upeta\in ]0,1[$, the set $D_{\upeta}$ 
is invariant under $T$}. Let us consider $z\in D_{\upeta}$.
From \eqref{eq:bound_for_tz0:prev},
we can deduce the following estimate
\begin{align}
\sum_{j=0}^{n-2}
|T^{(j)}z(t)| 
&\le \mathcal{R}(t)+
	\sum_{k=1}^n \Big(\|z\|_0\Big)^k
	\mathcal{L}_k(t)
\nonumber\\
&\le\mathcal{R}(t)+
	\upeta\mathcal{L}_1(t)+\upeta^2
	\left[\upeta\mathcal{L}_2(t)+
	\sum_{k=3}^n \Big(\upeta\Big)^{k-2}
	\mathcal{L}_k(t)\right]
\nonumber\\
&\le\mathcal{R}(t)+
	\upeta\mathcal{L}_1(t)+\upeta
\label{eq:estoamcion_estab}
\end{align}
in a right neighborhood of $\upeta=0$.
Now, by  \eqref{eq:integral_codition} we deduce that 
the first two terms on the  right side of \eqref{eq:estoamcion_estab}
tend to $0$ when $t\to\infty$. 
Hence, by \eqref{eq:estoamcion_estab} 
and (R2) we have that $\|Tz\|_0\le \upeta$,
or equivalently $Tz\in D_{\upeta}$ for all $z\in D_{\upeta}$.

\vskip 0.5cm
\noindent
{\bf (c)} {\it $T$ is a contraction for $\upeta\in ]0,1[$}. 
Let $z_1,z_2\in D_{\upeta}$. We can prove that
\begin{align*}
&\left|\sum_{|\boldsymbol{\alpha}|=1}^{n}
\Omega_{\boldsymbol{\alpha}}\Big(\mu,s,r_0(s),\ldots,r_{n-1}(s)\Big)
\Big[
\Big(z_1(s),\ldots,z_1^{(n-2)}(s)\Big)^{\boldsymbol{\alpha}}
-
\Big(z_2(s),\ldots,z_2^{(n-2)}(s)\Big)^{\boldsymbol{\alpha}}
\right|
\\
&
\quad
\le 
\sum_{|\boldsymbol{\alpha}|=1}^{n}
\Big|\Omega_{\boldsymbol{\alpha}}\Big(\mu,s,r_0(s),\ldots,r_{n-1}(s)\Big)\Big|
\left(|\boldsymbol{\alpha}|\upeta^{|\boldsymbol{\alpha}|-1}
	\sum_{i=0}^{n-2}\Big|z^{(i)}_1(s)-z^{(i)}_2(s)\Big|\right).
\end{align*}
Then, by the hypothesis  (R2)
and algebraic rearrangements,
we follow that
\begin{align*}
  &\sum_{j=0}^{n-2}\Big|T^{(j)} z_1(t)-T^{(j)} z_2(t)\Big|
   \\
   &
   \quad
   \le
   \int_{t_0}^{\infty}
   \sum_{j=0}^{n-2}\left|
	\frac{\partial^j g}{\partial t^j}(t,s) \right| 
	\sum_{|\boldsymbol{\alpha}|=1}^{n}
	\Big|\Omega_{\boldsymbol{\alpha}}\Big(\mu,s,r_0(s),\ldots,r_{n-1}(s)\Big)\Big|
	\left(
	|\boldsymbol{\alpha}|\upeta^{|\boldsymbol{\alpha}|-1}
	\sum_{i=0}^{n-2}|z^{(i)}_1(s)-z^{(i)}_2(s)|\right)
	ds
	\\
  &\quad
   \le
   \|z_1-z_2\|_0
   \int_{t_0}^{\infty}
   \sum_{j=0}^{n-2}\left|
	\frac{\partial^j g}{\partial t^j}(t,s) \right| 
	\sum_{|\boldsymbol{\alpha}|=1}^{n}
	|\boldsymbol{\alpha}|\upeta^{|\boldsymbol{\alpha}|-1}
	\Big|\Omega_{\boldsymbol{\alpha}}\Big(\mu,s,r_0(s),\ldots,r_{n-1}(s)\Big)\Big|
	ds
\\
  &\quad
   \le 	\|z_1-z_2\|_0
   \sum_{k=1}^{n}k\upeta^{k-1}\mathcal{L}_k(t)
   \quad\le\quad 	\|z_1-z_2\|_0\max\Big\{1,2\upeta,3\upeta^2,\ldots,n\upeta^{n-1}\Big\}
   \sum_{k=1}^{n}\mathcal{L}_k(t).
\end{align*}
Then, by application of \eqref{eq:integral_codition},
we deduce that $T$ is a contraction,
since, for an arbitrary $\upeta\in ]0,1[$, 
we have that $\max\Big\{1, 2\upeta, 3\upeta^2,\ldots,n\upeta^{n-1} \Big\}=1$.

Hence, from (a)-(c) and application of Banach fixed point theorem,
we deduce that there is a unique $z\in D_\upeta\subset
C_0^{n-2}([t_0,\infty[)$ solution of \eqref{eq:operator_equation}.

\subsection{Proof of Theorem~\ref{teo:asymptotic_ricatti}}
The proof is based on the operator equation~\eqref{eq:riccati_type_general}
and the invariant and contraction properties of $T$. Indeed,
let us first introduce some notation.
We denote by
 $z_{\mu}$ the solution of  the equation~\eqref{eq:riccati_type_general}
given by Theorem~\ref{teo:ricati}. Moreover, 
on $D_\upeta$ with $\upeta\in ]0,1/n[,$ we 
define recursively the sequence $\{\omega_{m}\}_{m\in\mathbb{N}}$  by
assuming that $\omega_{m+1}=T\omega_{m}$ with $\omega_{0}=0$. We note
that $\omega_{m}\to z_{\mu}$ when $m\to\infty$ as a consequence 
of the contraction property of $T$.

The Green
function $g$ defined on \eqref{eq:green_function} is given in terms
of $g_{\mu}$ and naturally,  the operator $T$
defined in \eqref{eq:operator_fix_point}
can be rewritten equivalently
as follows
\begin{align}
&Tz_{\mu}(t)
=\frac{1}{\Upsilon_0(\gamma_1,\ldots,\gamma_{n-1})}
\int_{t_0}^{\infty} g_{\mu}(t,s)
\mathbb{P}\Big(\mu,s,r_0(s),\ldots,r_{n-1}(s),
z_{\mu}(s),\ldots,z^{(n-2)}_{\mu}(s)\Big)ds
\nonumber\\
&
\;
=\frac{1}{\Upsilon_0(\gamma_1,\ldots,\gamma_{n-1})}
\nonumber\\
&
\;\;
\times \left\{
\begin{array}{l}
\displaystyle
\sum_{\ell=1}^{n-1} G_{\ell} 
\int_{t}^{\infty}
e^{-\gamma_\ell(t-s)}
\mathbb{P}\Big(\mu,s,r_0(s),\ldots,r_{n-1}(s),
z_{\mu}(s),\ldots,z^{(n-2)}_{\mu}(s)\Big)ds, 
\\
\hspace{8cm}\mathbb{G}_\mu\subset\mathbb{E}^{n-1}_1,
\\
\displaystyle
\sum_{\ell=1}^{k-1} G_{\ell}
\int_{t_0}^t
e^{-\gamma_\ell(t-s)} 
\mathbb{P}\Big(\mu,s,r_0(s),\ldots,r_{n-1}(s),
z_{\mu}(s),\ldots,z^{(n-2)}_{\mu}(s)\Big)ds 
\\
\displaystyle
\quad
+\sum_{\ell=k}^{n-1} G_{\ell}
\int_{t}^{\infty}
e^{-\gamma_\ell(t-s)}
\mathbb{P}\Big(\mu,s,r_0(s),\ldots,r_{n-1}(s),
z_{\mu}(s),\ldots,z^{(n-2)}_{\mu}(s)\Big)ds,
\\
\hspace{8cm}  \mathbb{G}_\mu\subset\mathbb{E}^{n-1}_k,\;\;\; k=\overline{2,n-1},
\\
\displaystyle
\sum_{\ell=1}^{n-1} G_{\ell}
\int_{t_0}^t
e^{-\gamma_\ell(t-s)}
\mathbb{P}\Big(\mu,s,r_0(s),\ldots,r_{n-1}(s),
z_{\mu}(s),\ldots,z^{(n-2)}_{\mu}(s)\Big)ds,
\\
\hspace{8cm}  \mathbb{G}_\mu\subset\mathbb{E}^{n-1}_n.
\end{array}
\right.
\label{eq:operator_fix_point_g1}
\end{align}
Here we have used the definition of Heaviside function, for instance in the case
$\mathbb{G}_\mu\subset\mathbb{E}^{n-1}_1$ the integration is on $[t,\infty[$ since  
$g_{\mu}(t,s)=0$ for $s\in [t_0,t]$.

Thus, the proof of \eqref{eq:asymptotic_teo:form} is reduced to get the
existence of the sequence $\{\Phi_m\}\subset\R_+$  such that
the following two assertions are valid:
\begin{align}
&
\sum_{j=0}^{n-2}|\omega^{(j)}_{m}(t)|
\leq
\Phi_m
\times
\left\{
\begin{array}{l}
\displaystyle
\int_{t}^{\infty}e^{-\beta(t-s)}
\Big|\Omega_{\boldsymbol{0}}
\Big(\mu,s,r_0(s),\ldots,r_{n-1}(s)\Big)\Big|ds,
\\
\hspace{4cm}\mathbb{G}_\mu\subset\mathbb{E}^{n-1}_1,
\quad \beta\in [\gamma_1,0[,
\\
\displaystyle
\int_{t_0}^{\infty}e^{-\beta(t-s)}
\Big|\Omega_{\boldsymbol{0}}\Big(\mu,s,r_0(s),
\ldots,r_{n-1}(s)\Big)\Big|ds,
\\
\hspace{3cm} 
\mathbb{G}_\mu\subset\mathbb{E}^{n-1}_k,\quad \beta\in [\gamma_k,0[,\quad 
 k=\overline{2,n-1},
\\
\displaystyle
\int_{t_0}^{t}e^{-\beta(t-s)}
\Big|\Omega_{\boldsymbol{0}}\Big(\mu,s,r_0(s),\ldots,r_{n-1}(s)\Big)\Big|ds,
\\
\hspace{3cm}
 \mathbb{G}_\mu\subset\mathbb{E}^{n-1}_n,\quad \beta\in ]0,\gamma_{n-1}],
\end{array}
\right.
\quad
\label{eq17}
\\
& \{\Phi_{m}\}
\mbox{ is convergent}.
\label{eq17:uniformly_bounded}
\end{align}
Hence,
to complete the proof of \eqref{eq:asymptotic_teo:form}, 
we proceed to prove \eqref{eq17} and \eqref{eq17:uniformly_bounded}.
Note that, the proof of \eqref{eq:asymptotic_teo:form} is 
concluded by passing to the limit the sequence $\{\Phi_{m}\}$ 
when $m\to\infty$ 
in the topology of $C^{n-2}_0([t_0,\infty]).$

\vspace{0.5cm}
\noindent
{\it  Proof of \eqref{eq17}}. 
We apply mathematical induction on $m$
and construct the sequences $\{\Phi_m\}$ for 
each $k=1,\ldots,n$. Indeed, the induction step
for $m=1$ is proved as follows.
From \eqref{eq:operator_fix_point_g1}, 
fact that $\omega_0=0,$ and (R1) we have that the 
following estimate 
\begin{align*}
&\sum_{j=0}^{n-2}|\omega^{(j)}_{1}(t)|   
=\sum_{j=0}^{n-2}|T \omega^{(j)}_0(t)|
\\
&
\le
\frac{1}{\Big|\Upsilon_0(\gamma_1,\ldots,\gamma_{n-1})\Big|}
\\
&
\;\;
\times \left\{
\begin{array}{l}
\displaystyle
\sum_{\ell=1}^{n-1} |G_{\ell}|(|\gamma_\ell|)^j
\int_{t}^{\infty}
e^{-\beta(t-s)}
\Big|\mathbb{P}\Big(\mu,s,r_0(s),\ldots,r_{n-1}(s),
0,\ldots,0\Big)\Big|ds, 
\\
\hspace{6cm}\mathbb{G}_\mu\subset\mathbb{E}^{n-1}_1,
\;\;\beta\in[\gamma_1,0[,
\\
\displaystyle
\sum_{\ell=1}^{k-1} |G_{\ell}|(|\gamma_\ell|)^j
\int_{t_0}^t
e^{-\beta(t-s)} 
\Big|\mathbb{P}\Big(\mu,s,r_0(s),\ldots,r_{n-1}(s),
0,\ldots,0\Big)\Big|ds 
\\
\displaystyle
\quad
+\sum_{\ell=k}^{n-1} |G_{\ell}|(|\gamma_\ell|)^j
\int_{t}^{\infty}
e^{-\beta(t-s)}
\Big|\mathbb{P}\Big(\mu,s,r_0(s),\ldots,r_{n-1}(s),
0,\ldots,0\Big)\Big|ds,
\\
\hspace{6cm}  \mathbb{G}_\mu\subset\mathbb{E}^{n-1}_k,
\;\;\beta\in[\gamma_k,0[,\;\;\; k=\overline{2,n-1},
\\
\displaystyle
\sum_{\ell=1}^{n-1} |G_{\ell}|(|\gamma_\ell|)^j
\int_{t_0}^t
e^{-\gamma_\ell(t-s)}
\Big|\mathbb{P}\Big(\mu,s,r_0(s),\ldots,r_{n-1}(s),
0,\ldots,0\Big)\Big|ds,
\\
\hspace{6cm}  \mathbb{G}_\mu\subset\mathbb{E}^{n-1}_n,\;\;\beta\in]0,\gamma_{n-1}[,
\end{array}
\right.
\end{align*}
is satisfied.
By application of   (R2) and the property 
\begin{align*}
 \mathbb{P}(\mu,s,r_0(s),\ldots,r_{n-1}(s),0,\ldots,0)
=\Omega_{\boldsymbol{0}}(\mu,s,r_0(s),\ldots,r_{n-1}(s)),
\end{align*}
we deduce that  \eqref{eq17} is valid with
$\Phi_{1}$ given on \eqref{eq:PHI_1_1}.
Then, the induction step for $m=1$ is valid.
Now, we prove the general induction step:
by assuming that \eqref{eq17} is valid for  
$m=h$, we will deduce that \eqref{eq17} is also valid for $m=h+1$.
Indeed, by the rewritten form of the operator $T$
given on \eqref{eq:operator_fix_point_g1}, we have that
\begin{align*}
&\sum_{j=0}^{n-2}|\omega^{(j)}_{h+1}(t)| 
=\sum_{j=0}^{n-2}|T\omega^{(j)}_{h}(t)|
\\ 
&\le 
\Phi_1\times
\left\{
\begin{array}{l}
\displaystyle
\int_{t}^{\infty}
e^{-\gamma_1(t-s)}
\Big|\mathbb{P}\Big(\mu,s,r_0(s),\ldots,r_{n-1}(s),
\omega_{h}(s),\ldots,\omega_{h}^{(n-2)}(s)\Big)\Big|ds, 
\quad \mathbb{G}_\mu\subset\mathbb{E}^{n-1}_1,
\\
\displaystyle
\int_{t_0}^t
e^{-\gamma_\ell(t-s)} 
\mathbb{P}\Big(\mu,s,r_0(s),\ldots,r_{n-1}(s),
\omega_{h}(s),\ldots,\omega_{h}^{(n-2)}(s)\Big)ds 
\\
\displaystyle
\qquad\quad
+
\int_{t}^{\infty}
e^{-\gamma_\ell(t-s)}
\mathbb{P}\Big(\mu,s,r_0(s),\ldots,r_{n-1}(s),\omega_{h}(s),\ldots,\omega_{h}^{(n-2)}(s)\Big)ds,
\\
\hspace{9.0cm}\mathbb{G}_\mu\subset\mathbb{E}^{n-1}_k,  k=\overline{2,n-1},
\\
\displaystyle
\int_{t_0}^t
e^{-\gamma_n(t-s)}
\Big|\mathbb{P}\Big(\mu,s,r_0(s),\ldots,r_{n-1}(s),\omega_{h}(s),\ldots,\omega_{h}^{(n-2)}(s)\Big)\Big|ds.
\quad \mathbb{G}_\mu\subset\mathbb{E}^{n-1}_n,
\end{array}
\right.
\end{align*}
with $\Phi_1$ defined on \eqref{eq:PHI_1_1}.
We study separated the three cases: 
$\mathbb{G}_\mu\subset\mathbb{E}^{n-1}_1,$ 
$\mathbb{G}_\mu\subset\mathbb{E}^{n-1}_k$
for $k=\overline{2,n-1},$ and $\mathbb{G}_\mu\subset\mathbb{E}^{n-1}_n$, respectively.
Firstly, for $\mathbb{G}_\mu\subset\mathbb{E}^{n-1}_1$, 
by application of the hypothesis (R3) and $\beta\in[\gamma_1,0[$, we have that
\begin{align}
&\sum_{j=0}^{n-2}|\omega^{(j)}_{h+1}(t)|
\le\Phi_{1}
\int_{t}^{\infty} 
\exp(-\gamma_1(t-s)) 
\Big|\mathbb{P}\Big(\mu,s,r_0(s),\ldots,r_{n-1}(s),
\omega_{h}(s),\ldots,w^{(n-2)}_{h}(s)\Big)\Big|ds
\nonumber\\
&
\quad
\le\Phi_{1}
\int_{t}^{\infty} \exp(-\gamma_1(t-s)) 
\sum_{|\boldsymbol{\alpha}|=0}^n
\Big|\Omega_{\boldsymbol{\alpha}}\Big(\mu,s,r_0(s),\ldots,r_{n-1}(s)\Big)\Big|
\Big|\Big(\omega_h(s),\ldots,\omega_h^{(n-2)}(s)\Big)^{\boldsymbol{\alpha}}\Big|ds
\nonumber\\
&
\quad
\le
\Phi_{1}
\Bigg\{\int_{t}^{\infty} \exp(-\gamma_1(t-s))
\Big|\Omega_{\boldsymbol{0}}\Big(\mu,s,r_0(s),\ldots,r_{n-1}(s)\Big)\Big|ds
\nonumber\\
&
\hspace{1cm}
+
\int_{t}^{\infty} \exp(-\gamma_1(t-s))
\sum_{|\boldsymbol{\alpha}|=1}^n
\Big|\Omega_{\boldsymbol{\alpha}}\Big(\mu,s,r_0(s),\ldots,r_{n-1}(s)\Big)\Big|
\Big(\upeta^{|\boldsymbol{\alpha}|-1}
\sum_{j=0}^{n-2}|\omega^{(j)}_{h}(s)|
\Big)ds
\Bigg\}
\nonumber\\
&
\quad
\le
\Phi_{1}
\Bigg\{\int_{t}^{\infty} \exp(-\gamma_1(t-s))
\Big|\Omega_{\boldsymbol{0}}\Big(\mu,s,r_0(s),\ldots,r_{n-1}(s)\Big)\Big|ds
+\Phi_{h}\int_{t}^{\infty} \exp(-\gamma_1(t-s))
\nonumber\\
&
\hspace{0.5cm}
\times
\sum_{|\boldsymbol{\alpha}|=1}^n
\Big|\Omega_{\boldsymbol{\alpha}}\Big(\mu,s,r_0(s),\ldots,r_{n-1}(s)\Big)\Big|
\int_{s}^{\infty}\!\!\! \exp(-\beta(t-\tau))
\Big|\Omega_{\boldsymbol{0}}\Big(\mu,\tau,r_0(\tau),\ldots,r_{n-1}(\tau)\Big)\Big|d\tau
ds\Bigg\}
\nonumber\\
&
\quad
\le
\Phi_{1}
\Bigg\{1+\Phi_{h}\int_{t}^{\infty} \exp(-\gamma_1(t-s))
\sum_{|\boldsymbol{\alpha}|=1}^n
\Big|\Omega_{\boldsymbol{\alpha}}\Big(\mu,s,r_0(s),\ldots,r_{n-1}(s)\Big)\Big|
ds\Bigg\}
\nonumber\\
&
\hspace{1.5cm}
\times
\int_{t}^{\infty} \exp(-\beta(t-\tau))
\Big|\Omega_{\boldsymbol{0}}\Big(\mu,\tau,r_0(\tau),\ldots,r_{n-1}(\tau)\Big)\Big|d\tau
\nonumber\\
&
\quad
\le
\Phi_{1}
\Big(1+\sigma_{\gamma_1}\Phi_{h}\Big)
\int_{t}^{\infty} \exp(-\beta(t-s))
\Big|\Omega_{\boldsymbol{0}}\Big(\mu,s,r_0(s),\ldots,r_{n-1}(s)\Big)\Big|ds.
\label{eq:induccionk1}
\end{align}
For $\mathbb{G}_\mu\subset\mathbb{E}^{n-1}_k$
with $k=\overline{2,n-1},$ 
from (H3) and $\beta\in[\gamma_k,0[$,
by a similar arguments 
we deduce the estimate 
\begin{align}
\sum_{j=0}^{n-2}|\omega^{(j)}_{h+1}(t)|
 &\le
\Phi_{1}\left[
\int_{t_0}^t
e^{-\gamma_\ell(t-s)} 
\mathbb{P}\Big(\mu,s,r_0(s),\ldots,r_{n-1}(s),\omega_h(s),\ldots,\omega^{(n-2)}_h(s)\Big)ds 
\right.
\nonumber\\
&\displaystyle
\qquad\quad
\left.
+
\int_{t}^{\infty}
e^{-\gamma_\ell(t-s)}
\mathbb{P}\Big(\mu,s,r_0(s),\ldots,r_{n-1}(s),\omega_h(s),\ldots,\omega_h^{(n-2)}(s)\Big)ds
\right]
\nonumber\\
&
\quad
\le
\Phi_{1}
\Bigg\{1+\Phi_{h}\int_{t_0}^{\infty} \exp(-\gamma_k(t-s))
\sum_{|\boldsymbol{\alpha}|=1}^n
\Big|\Omega_{\boldsymbol{\alpha}}\Big(\mu,s,r_0(s),\ldots,r_{n-1}(s)\Big)\Big|
ds\Bigg\}
\nonumber\\
&
\hspace{1.5cm}
\times
\int_{t_0}^{\infty} \exp(-\beta(t-\tau))
\Big|\Omega_{\boldsymbol{0}}\Big(\mu,\tau,r_0(\tau),\ldots,r_{n-1}(\tau)\Big)\Big|d\tau
\nonumber\\
&
\quad
\le
\Phi_{1}
\Big(1+\sigma_{\gamma_k}\Phi_{h}\Big)
\int_{t_0}^{\infty} \exp(-\beta(t-s))
\Big|\Omega_{\boldsymbol{0}}\Big(\mu,s,r_0(s),\ldots,r_{n-1}(s)\Big)\Big|ds.
\label{eq:induccionkmedio}
\end{align}
The case $\mathbb{G}_\mu\subset\mathbb{E}^{n-1}_n$ is similar to the case 
$\mathbb{G}_\mu\subset\mathbb{E}^{n-1}_1$ and we get
\begin{align}
\sum_{j=0}^{n-2}|\omega^{(j)}_{h+1}(t)|
\le
\Phi_{1}
\Big(1+\sigma_{\gamma_n}\Phi_{h}\Big)
\int_{t_0}^t \exp(-\beta(t-s))
\Big|\Omega_{\boldsymbol{0}}\Big(\mu,s,r_0(s),\ldots,r_{n-1}(s)\Big)\Big|ds,
\label{eq:induccionkn}
\end{align}
for $\beta\in]0,\gamma_{n-1}[.$
Then, from \eqref{eq:induccionk1}-\eqref{eq:induccionkn},
we deduce that the induction process is valid and 
\eqref{eq17} is satisfied with
\begin{align}
\Phi_{h}=   \Phi_{1}\Big(1+\Phi_{h-1}\;\sigma_{\gamma_k}\Big),
\quad k=1,\ldots,n,
\label{eq:seq_pimi}
\end{align}
where  $\Phi_{1}$ is defined on \eqref{eq:PHI_1_1}.

\vspace{0.5cm}
\noindent
{\it  Proof of \eqref{eq17:uniformly_bounded}}.
From \eqref{eq:seq_pimi},
using recursively the definition of $\Phi_{h-1},\ldots,\Phi_{2}$,
we can rewrite $\Phi_{h}$ as the sum of the terms of 
a geometric progression where the common ratio is given
by $\sigma_{\gamma_i}\Phi_{1}$. Then, 
the hypothesis (R3) implies  that
$\sigma_{\gamma_i}\Phi_{1}\in ]0,1[$,
and we can deduce that         
\begin{align*}
\lim_{m\to\infty}
\Phi_{m}&=\Phi_{1}\lim_{m\to\infty}\sum_{i=0}^{m-1}
	\Big(\Phi_{1}\sigma_{\gamma_i}\Big)^i
        =\Phi_{1}
        \lim_{m\to\infty}
        \frac{\Big[(\Phi_{1}\sigma_{\gamma_i})^m-1\Big]}
        {\Phi_{1}\sigma_{\gamma_i}-1}
        =
        \frac{\Phi_{1}}{1-\Phi_{1}\sigma_{\gamma_i}}
        =\Phi_{\gamma_i}>0
\end{align*}
or equivalently $\Phi_{m}$ converges to $\Phi_{\gamma_i}.$

\subsection{Proof of Theorem~\ref{teo:changeofvariable}}

[{\it (a)}] In order to prove \eqref{eq:ricati_enesima}, we firstly construct
a formula for 
$y^{(j)}$ and then we use that formula in \eqref{eq:intro_uno_intro}.
Indeed, we have that the $j-$th derivative for  $y$ is given by
\begin{align}
y^{(j)}(t)=
\left\{
\begin{array}{ll}
 (z(t)+\mu)y(t), & j=1,
 \\
 \Big[(z(t)+\mu)^{(1)}+(z(t)+\mu)^2\Big]y(t), & j=2,
 \\
 \Big[\sum_{m=0}^{j-3}S_{m,j}(t)
 + z^{j-2}(t)z^{(1)}(t)+(z(t)+\mu)^{j}\Big]y(t), &
 j\ge 3.
\end{array}
\right.
\label{eq:leibiniz:j} 
\end{align}
The proof of \eqref{eq:leibiniz:j} for $j=1,2$ is given by
direct differentiation of $y$ in \eqref{eq:general_change_var}.
Now, for $j\ge 3$, we proceed inductively by using  the Leibniz formula for 
differentiation. Indeed, by using the relation for $y^{(1)}$ we have that
\begin{align}
y^{(j)}(t)&=(y^{(1)}(t))^{(j-1)}=\Big((z(t)+\mu)y(t)\Big)^{(j-1)}
=\sum_{\ell_1=0}^{j-1}{j-1\choose \ell_1}(z(t)+\mu)^{(\ell_1)} y^{(j-1-\ell_1)}(t).
\nonumber\\
&=
\Big[(z(t)+\mu)^{(j-1)} +(j-1)(z(t)+\mu)^{(j-2)} (z(t)+\mu)\Big]y(t)
\nonumber\\
&
\quad+
\sum_{\ell_1=0}^{j-3}{j-1\choose \ell_1}(z(t)+\mu)^{(\ell_1)} y^{(j-\ell_1-1)}(t).
\label{eq:leibiniz:1}
\end{align}
The order of derivatives for  $y$ is strictly decreasing wit respect
to $\ell_1$, since $j-0-1>j-1-1>\ldots>j-(j-3)-1=2$. Then, if $j-1>2 $ we can apply again 
the Leibniz formula to compute $ y^{(j-\ell_1-1)}=(y')^{(j-\ell_1-2)}$ and 
from \eqref{eq:leibiniz:1} we get
\begin{align}
&y^{(j)}(t)=
\Big[(z(t)+\mu)^{(j-1)} +(j-1)(z(t)+\mu)^{(j-2)} (z(t)+\mu)\Big]y(t)
\nonumber\\
&
\quad
+
\sum_{\ell_1=0}^{j-3}{j-1\choose \ell_1}
(z(t)+\mu)^{(\ell_1)}
\Big[(z(t)+\mu)^{(j-\ell_1-2)} +(j-\ell_1-2)(z(t)+\mu)^{(j-\ell_1-3)} (z(t)+\mu)\Big]y(t)
\nonumber\\
&
\quad
+\sum_{\ell_1=0}^{j-4}
\sum_{\ell_2=0}^{j-\ell_1-4}
{j-1\choose \ell_1}
{j-\ell_1-2\choose \ell_2}
(z(t)+\mu)^{(\ell_1)}
(z(t)+\mu)^{(\ell_2)}
y^{(j-\ell_1-\ell_2-2)}(t).
\label{eq:leibiniz:3}
\end{align}
Similarly, observing the order of derivatives for $y$,
we deduce that
if $j-2>2$ we can use  the Leibniz formula  to compute $y^{(j-\ell_1-\ell_2-2)}$.
We note that, we can apply this strategy $j-2$ times to obtain the desired result 
given in \eqref{eq:leibiniz:j} for $j\ge 3$ and conclude the proof
of item (a).

Replacing  \eqref{eq:leibiniz:j}
in \eqref{eq:intro_uno_intro} we  deduce that
\begin{align*}
&y^{(n)}(t)+\sum_{i=0}^{n-1}[a_i+r_i (t)] y^{(i)}(t)=
y^{(n)}(t)+\sum_{i=3}^{n-1}[a_i+r_i (t)] y^{(i)}(t)+\sum_{i=0}^{2}[a_i+r_i (t)] y^{(i)}(t)
\\
&\;
=\Bigg\{
\Big[z^{(n-1)}(t)+n\mu^{(n-1)}z(t)+(n-1)\mu\;z^{(n-2)}(t)\Big]+\mu^{n}
+\Big[(n-1)z^{(n-2)}(t)z(t)+\sum_{m=1}^{n-3}S_{m,n}(t)
\\
&\qquad
+z^{(n-1)}(t)z^{(1)}(t)+\sum_{j=0}^{n-2}{n\choose j}
z^{n-j}(t)\mu^{j}\Big]
\\
&\quad
+
\Big[\sum_{i=3}^{n-1} a_i\Big(z^{(i-1)}(t)+(i-1)z^{(i-2)}(t)\mu+i\mu^{i-1}z(t)\Big)\Big]
+\sum_{i=3}^{n-1} a_i\mu^i
+\Big[
\sum_{i=3}^{n-1} a_i\Big((i-1)z^{(i-2)}(t)z(t)
\\
&\qquad\quad
+\sum_{m=1}^{i-3}S_{m,i}(t)+z^{(i-1)}(t)z^{(1)}(t)+\sum_{j=1}^{i-2}{i\choose j}
z^{i-j}(t)\mu^{j}\Big) 
+
r_i(t)\Big(\sum_{m=0}^{i-3}S_{m,i}(t)+z^{(i-1)}(t)z^{(1)}(t)
\\
&\hspace{2cm}
+(z(t)+\mu)^{i}\Big)\Big]
+
\Big[a_2z^{(1)}(t)+2a_2 \mu z(t)+a_1 z(t)\Big]
+\Big[a_2\mu^2+a_1\mu+a_0\Big]
\\
&\quad
+\Big[a_2z^2(t)+r_2(t)\Big(z^{(1)}(t)+(z(t)+\mu)^2\Big)
+r_1(t)\Big(z(t)+\mu\Big)+r_0(t)\Big]
\Bigg\}y(t)=0.
\end{align*}
Thus, using the fact that $\mu$ is a characteristic root we have that
$z$ is a solution of \eqref{eq:ricati_enesima}.

\vspace{0.5cm}
\noindent
[{\it (b)}] The proof of \eqref{eq:charset_ricca} is a consequence of the following
claim: If
$\lambda_i,\lambda_j\in \Lambda$ for $i\not=j$,
then $\lambda_j-\lambda_i$ is a 
root of the characteristic polynomial associated to
the differential equation \eqref{eq:ricati_enesima}
when $\mu=\lambda_i$ and the right hand side is zero.
Indeed, using the identity 
\begin{align*}
\sum_{k=0}^{i-1}{i\choose k}(u-v)^{i-1-k}v^k=\sum_{k=0}^{i-1}u^{i-1-k}v^k
\end{align*}
we deduce that
\begin{align*}
&(u^n+\sum_{i=0}^{n-1}a_iu^i)-(v^n+\sum_{i=0}^{n-1}a_iv^i)
=(u-v)\Big[\sum_{k=0}^{n-1}u^{n-1-k}v^k
+\sum_{i=2}^{n-1}a_i \sum_{k=0}^{i-1}u^{i-1-k}v^k +a_1
\Big]
\\
&
\quad
=(u-v)\Big[\sum_{k=0}^{n-1}{n\choose k}(u-v)^{n-1-k}v^k
+\sum_{i=2}^{n-1} \sum_{k=0}^{i-1}{i\choose k}a_i(u-v)^{i-1-k}v^k+a_1
\Big]
\\
&
\quad
=(u-v)\Big[
(u-v)^{(n-1)}+\sum_{i=2}^{n-1}\Big[\sum_{k=i}^{n}a_k {k \choose i} v ^{k-i} \Big]
(u-v)^{(i-1)}
+\Big(n\mu ^{n-1}+\sum_{i=1}^{n-1}i a_i \mu ^{i-1}\Big) (u-v)
\Big]
\end{align*}
and by selecting $u=\lambda_j$ and $v=\lambda_i$ we  can prove the claim.

\vspace{0.5cm}
\noindent
[{\it (c)}] We note that the equation \eqref{eq:ricati_enesima} is of the type
\eqref{eq:riccati_type_general} where (R1) and (R2) are satisfied, since
from the item (b) of Theorem~\ref{teo:changeofvariable} we have that (H1)  implies (R1)
and also clearly (H2) implies (R2).
Thus, from Theorem~\ref{teo:ricati}, we deduce the existence and uniqueness
of $z_{\lambda_i}$ belong to $C_0^{n-2}([t_0,\infty[)$ 
satisfying \eqref{eq:ricati_enesima}. 

\vspace{0.5cm}
\noindent
[{\it (d)}]  We note that $\sum_{|\boldsymbol{\alpha}|=1}^n
\Big|\Omega_{\boldsymbol{\alpha}}\Big(\mu,s,r_0(s),\ldots,r_{n-1}(s)\Big)\Big|$
is the sum on the coefficients of $\mathbb{P}$ except of the independent term
and $\mathbb{H}(s,\mu)$ is the of the coefficients
of $\mathbb{F}$ except of the independent term. Also observe that 
$\mathbb{H}(s,\mu)$ can be obtained from \eqref{eq:function:FF}
by considering that $z=z'=\ldots=z^{n-2}=1$. Then, we have that
(H3) implies (R3).

\subsection{Proof of Theorem~\ref{teo:main}}
From Theorem~\ref{teo:changeofvariable}, 
we have that the fundamental system of solutions for
\eqref{eq:intro_uno_intro} is given by \eqref{eq:fundam_sist_sol}.
Moreover, from \eqref{eq:leibiniz:j} we deduce that
\begin{align}
\frac{y^{(j)}_i(t)}{y_i(t)}=
\left\{
\begin{array}{ll}
 (z_{\lambda_i}(t)+\lambda_i), & j=1,
 \\
 \Big[(z_{\lambda_i}(t)+\lambda_i)^{(1)}+(z_{\lambda_i}(t)+\lambda_i)^2\Big], & j=2,
 \\
 \Big[\sum_{m=0}^{j-3}S_{m,j}(t)
 + z_{\lambda_i}^{j-2}(t)z_{\lambda_i}^{(1)}(t)+(z_{\lambda_i}(t)+\lambda_i)^{j}\Big], &
 j\ge 3,
\end{array}
\right.
\label{eq:leibiniz:dos} 
\end{align}
Now, using the fact that $z_{\lambda_i}\in C_0^2([t_0,\infty[)$
is a solution of \eqref{eq:fundam_sist_sol} with
$\mu=\lambda_i$, we deduce
the proof of \eqref{eq:asymptotic_perturbed_w0}.
By the definition of the $W[y_1,\ldots,y_n]$,
some algebraic rearrangements and 
\eqref{eq:asymptotic_perturbed_w0}, we deduce \eqref{eq:asymptotic_perturbed_w}.

In order to prove \eqref{eq:asymptotic_perturbed_3_4} we 
use the relations \eqref{eq:fundam_sist_sol}, \eqref{eq:operator_fix_point}, 
\eqref{eq:operator_fix_point_g1}
and the identity
\begin{align}
&\int_{t_0}^t e^{-a\tau}\int_{\tau}^\infty e^{as} H(s)dsd\tau
\nonumber\\
&
\hspace{2cm}
=-\frac{1}{a}\left[
\int_{t}^\infty e^{-a(t-s)}H(s)ds-\int_{t_0}^\infty e^{-a(t_0-s)}H(s)ds
\right]
+\frac{1}{a}\int_{t_0}^t H(\tau)d\tau.
\label{eq:idenaux_pro}
\end{align}
Indeed, the relation \eqref{eq:asymptotic_perturbed_3_4}
implies that
\begin{align}
y_i(t)=\exp\Big(\int_{t_0}^t(\lambda_i+z_i(\tau))d\tau\Big)
=e^{\lambda_i(t-t_0)}\exp\Big(\int_{t_0}^t z_{\lambda_i}(\tau)d\tau\Big).
\label{eq:solfun_paso1}
\end{align}
Now, from (H1) and Theorem~\ref{teo:changeofvariable}-(b) we can deduce that
the sets defined on \eqref{eq:notation:G} and \eqref{eq:notation:E}
satisfy the inclusion $\mathbb{G}_{\lambda_i}\subset \mathbb{E}^{n-1}_i$.
Then, we can rewrite \eqref{eq:operator_fix_point_g1} with $\mathbb{F}$ instead of $\mathbb{P}$
and $\mu=\lambda_i$, i.e.
\begin{align*}
&Tz_{\lambda_i}(t)
=\frac{1}{\Upsilon_0\Big(\gamma_1(\lambda_i),\ldots,\gamma_{n-1}(\lambda_i)\Big)}
\nonumber\\
&
\;\;
\times \left\{
\begin{array}{l}
\displaystyle
\sum_{\ell=2}^{n} G_{\ell} 
\int_{t}^{\infty}
e^{-(\lambda_\ell-\lambda_1)(t-s)}
\mathbb{F}\Big(\lambda_1,s,r_0(s),\ldots,r_{n-1}(s),
z_{\lambda_1}(s),\ldots,z^{(n-2)}_{\lambda_1}(s)\Big)ds, \;\; i=1
\\
\displaystyle
\sum_{\ell=1}^{i-1} G_{\ell}
\int_{t_0}^t
e^{-(\lambda_\ell-\lambda_i)(t-s)} 
\mathbb{F}\Big(\lambda_i,s,r_0(s),\ldots,r_{n-1}(s),
z_{\lambda_i}(s),\ldots,z^{(n-2)}_{\lambda_i}(s)\Big)ds 
\\
\displaystyle
\quad
+\sum_{\ell=i}^{n-1} G_{\ell}
\int_{t}^{\infty}
e^{-(\lambda_{\ell+1}-\lambda_i)(t-s)}
\mathbb{F}\Big(\lambda_i,s,r_0(s),\ldots,r_{n-1}(s),
z_{\lambda_i}(s),\ldots,z^{(n-2)}_{\lambda_i}(s)\Big)ds,
\\
\hspace{10cm} i=\overline{2,n-1},
\\
\displaystyle
\sum_{\ell=1}^{n-1} G_{\ell}
\int_{t_0}^t
e^{-(\lambda_{\ell}-\lambda_n)(t-s)}
\mathbb{F}\Big(\lambda_n,s,r_0(s),\ldots,r_{n-1}(s),
z_{\lambda_n}(s),\ldots,z^{(n-2)}_{\lambda_n}(s)\Big)ds,
\;\; i=n,
\end{array}
\right.
\end{align*}
where 
$\displaystyle G_\ell=(-1)^\ell\Upsilon_{\ell}
\Big(\gamma_1(\lambda_i),\ldots,\gamma_{n-1}(\lambda_i)\Big)$
and 
$\mathbb{F}$ is defined on \eqref{eq:function:FF}. Using
\eqref{eq:operator_fix_point}, we have that
\begin{align*}
&\int_{t_0}^t z_{\lambda_i}(\tau)d\tau
=\int_{t_0}^t Tz_{\lambda_i}(\tau)d\tau
=\frac{1}{\Upsilon_0(\lambda_1,\ldots,\lambda_n)}
\\
&
\;\;
\times \left\{
\begin{array}{l}
\displaystyle
\sum_{\ell=2}^{n} G_{\ell} 
\int_{t_0}^t\int_{\tau}^{\infty}
e^{-(\lambda_\ell-\lambda_1)(\tau-s)}
\mathbb{F}\Big(\lambda_1,s,r_0(s),\ldots,r_{n-1}(s),
z_{\lambda_1}(s),\ldots,z^{(n-2)}_{\lambda_1}(s)\Big)dsd\tau, \;\; i=1
\\
\displaystyle
\sum_{\ell=1}^{i-1} G_{\ell}
\int_{t_0}^t\int_{t_0}^\tau
e^{-(\lambda_\ell-\lambda_i)(\tau-s)} 
\mathbb{F}\Big(\lambda_i,s,r_0(s),\ldots,r_{n-1}(s),
z_{\lambda_i}(s),\ldots,z^{(n-2)}_{\lambda_i}(s)\Big)dsd\tau 
\\
\displaystyle
\quad
+\sum_{\ell=i}^{n-1} G_{\ell}
\int_{t_0}^t\int_{\tau}^{\infty}
e^{-(\lambda_{\ell+1}-\lambda_i)(\tau-s)}
\mathbb{F}\Big(\lambda_i,s,r_0(s),\ldots,r_{n-1}(s),
z_{\lambda_i}(s),\ldots,z^{(n-2)}_{\lambda_i}(s)\Big)dsd\tau,
\\
\hspace{10cm} i=\overline{2,n-1},
\\
\displaystyle
\sum_{\ell=1}^{n-1} G_{\ell}
\int_{t_0}^t\int_{t_0}^\tau
e^{-(\lambda_{\ell}-\lambda_n)(\tau-s)}
\mathbb{F}\Big(\lambda_n,s,r_0(s),\ldots,r_{n-1}(s),
z_{\lambda_n}(s),\ldots,z^{(n-2)}_{\lambda_n}(s)\Big)dsd\tau,
\;\; i=n,
\end{array}
\right.
\\
&=
\Big[\displaystyle\prod_{1\le j\le n, j\not=i}(\lambda_j-\lambda_i)\Big]^{-1}
	\int_{t_0}^t
	\mathbb{F}\Big(\lambda_i,s,r_0(s),\ldots,r_{n-1}(s),
z_{\lambda_i}(s),\ldots,z^{(n-2)}_{\lambda_i}(s)\Big)ds
	+o(1).
\end{align*}
Then, replacing in  \eqref{eq:solfun_paso1} we have that 
\eqref{eq:asymptotic_perturbed_3_4} is valid.

\section*{Acknowledgments}

The authors would like to thank the support of 
research projects DIUBB  172409 GI/C 
and DIUBB  183309~4/R at
Universidad del B{\'\i}o-B{\'\i}o, Chile.


\begin{thebibliography}{99}

%
%
%

\bibitem{bellman_book} {\textsc{R. Bellman},
 \textit{Stability Theory of Differential Equations},
 McGraw-Hill Book Company, Inc., New York-Toronto-London, 1953.}

%

\bibitem{coronel_huancas_pinto} \textsc{A. Coronel},
  \textit{Asymptotic integration of a linear fourth order differential equation of Poincar\'e type.}
 Electron. J. Qual. Theory Differ. Equ.  2015, No. 76, 24 pp.
%
%
%
%
%
%
%
%
%
\bibitem{eastham_book} {\textsc{M.S.P. Eastham},
 \textit{ The Asymptotic Solution of Linear Differential Dystems,
 Applications of the Levinson theorem}.
  London Mathematical Society Monographs, volume~4,
  Oxford University Press, New York, 1989.}
%
%

\bibitem{figueroa_2006} {\textsc{P. Figueroa, M. Pinto},
 Asymptotic expansion of the variable eigenvalue associated to second
  order differential equations, \textit {Nonlinear Stud.}, {\bf 13(3)} (2006),
  261--272.}

\bibitem{figueroa_2008} {\textsc{P. Figueroa, M. Pinto} { Riccati
equations and nonoscillatory solutions of third order
  differential equations}, \textit{ Dynam. Systems Appl.}, {\bf 17(3-4)} (2008),
  459--475.}
%
%
%
%
%
\bibitem{harris_lutz_1977} {\textsc{W.A. Harris Jr., D.A. Lutz}, {
A unified theory of asymptotic integration},
  \textit{J. Math. Anal. Appl.}, \textbf { 57(3)} (1977),
  571--586.}


\bibitem{hartman_1948} {\textsc{P. Hartman}, { Unrestricted
solution fields of almost-separable differential
  equations},
\textit{ Trans. Amer. Math. Soc.},  \textbf { 63} (1948),
560--580.}

\bibitem{hartman_wintner_1955} {\textsc{P. Hartman, A. Wintner}, {
Asymptotic integrations of linear differential equations},
 \textit{ American Journal of Mathematics}, \textbf{77(1)} (1955),
 45--86.}

%
%

\bibitem{levinson1948} {\textsc{N. Levinson}, { The asymptotic
nature of solutions of linear systems of differential
  equations},
 \textit{ Duke Mathematical Journal}, \textbf { 15(1)} (1948),
 111--126.}


\bibitem{perron_1909} {\textsc{O. Perron},{ Ber einen satz des henr
Poincar\'e},
 \textit{J. Reine Angew. Math.}, \textbf {136} (1909), 17--37.}

%
%
%
%

\bibitem{poincare_1885} {\textsc{H. Poincar\'e},{
 Sur les equations lineaires aux differentielles ordinaires et
  aux differences finies},
 \textit{  Amer. J. Math.}, \textbf { 7(3)} (1885), 203--258.}

%
%
%
%
%
%
%
%

\end{thebibliography}
\end{document}